\documentclass{svproc}
\usepackage{url}

\usepackage{amsmath}
\usepackage{amssymb}
\usepackage{amsfonts}
\usepackage{latexsym}
\usepackage{mathtools}
\usepackage{empheq}
\usepackage{color}

\newcommand{\capE}{\mathcal{E}}

\newcommand{\mathN}{\mathbb{N}}

\newcommand{\mathR}{\mathbb{R}}

\begin{document}
\mainmatter              
\title{Local and Nonlocal Liquid Drop Models}
\titlerunning{Liquid Drop Model}  
%
\author{Matteo Novaga\inst{1} \and Fumihiko Onoue\inst{2}
}
\authorrunning{M. Novaga and F. Onoue} 
%
%
\institute{Universit\`a di Pisa, Largo Bruno Pontecorvo 5, 56127 Pisa, Italy\\
\email{matteo.novaga@unipi.it}
\and
Technische Universit\"at M\"unchen, Boltzmannstrasse 3, 85748 Garching, Germany\\
\email{fumihiko.onoue@tum.de}}

\maketitle              

\begin{abstract}
We consider some extensions of Gamow's liquid drop model for an atomic nucleus. We present a review of the classical  model and then
we illustrate some recent developments on a nonlocal variant, where the perimeter term is replaced by the fractional perimeter.

\end{abstract}

\section{Introduction}
We review some recent developments on the liquid drop model introduced by George Gamow in \cite{Gamow} in 1930s, to explain the behavior of atomic nuclei and predict the phenomena of nuclear fission.

In this model, the attractive short-range nuclear force produces a surface tension due to lower nucleon density near the nucleus boundary. Meanwhile, the presence of protons, which is positively charged, produces the repulsive Coulomb force. Treating the collection of protons and neutrons inside an atomic nucleus as an incompressible uniformly charged fluid, the model can be written mathematically as the following energy:
\begin{equation*}
	\capE(E) \coloneqq |\partial E| + \int_{E}\int_{E} \frac{1}{|x-y|}\,dx\,dy
\end{equation*}
where the nucleus $E \subset \mathR^3$ is a smooth set with fixed volume $|E| = m$ and $|\partial E|$ is the area of the boundary of $E$. The volume $m$ is a parameter proportional to the number of neutrons in a nucleus. Then the ground state of a nucleus with a given number of nucleons is the minimizer of $\capE$, namely, the set $E$ that attains the least energy,
\begin{equation*}
	E[m] \coloneqq \inf\left\{ \capE(E) \coloneqq |\partial E| +  \int_{E}\int_{E}\frac{1}{|x-y|}\,dx\,dy \mid E \subset \mathR^3, \, |E| = m \right\},
\end{equation*}
for a given volume $m>0$. 

From a mathematical point of view, Problem $E[m]$ for $m > 0$ can be generalized into the following minimization problem:
\begin{equation}\label{classicalLiquidDropModel}
	\inf\left\{ \capE_{g}(E) \coloneqq P(E) +  V_g(E) \mid |E| = m \right\}
\end{equation}
where $P(E)$ is the classical perimeter of $E \subset \mathR^N$ in the sense of De Giorgi, and $V_g(E)$ is the generalized Riesz potential term of $E \subset \mathR^N$ given as
\begin{equation}\label{definitionGeneralRieszPotential}
	V_g(E) \coloneqq \int_{E}\int_{E} g(x-y) \,dx\,dy.
\end{equation}
A typical example of $g$ is given by $|x|^{-(N-\alpha)}$ for $\alpha \in (0,\,N)$ and one can easily notice that, if $N=3$ and $\alpha = N-1$, then Problem \eqref{classicalLiquidDropModel} is equivalent to Problem $E[m]$ for $m>0$. From the isoperimetric inequality of the classical perimeter, one sees that the ball is the unique minimizer of the classical perimeter among sets with fixed volume. In contrast, by the Riesz arrangement inequality, one sees that the ball is the maximizer of the Riesz potential. Hence, in Problem \eqref{classicalLiquidDropModel}, the non-trivial competition between the classical perimeter and Riesz potential term occurs. Moreover, when the kernel $g$ is given by $|x|^{-(N-\alpha)}$ with $\alpha \in (0,\,N)$, one can heuristically observe the existence of minimizers for small volumes and the non-existence of minimizers for large volumes. Indeed, by rescaling $E \mapsto F \coloneqq \lambda^{-1} \, E$ with $\lambda^N \coloneqq |B_1|^{-1} |E|$, we have that 
\begin{equation*}
	\capE_{g}(E) = \lambda^{N-1} \, P(F) + \lambda^{N+\alpha} \, V_{\alpha}(F) = \lambda^{N-1} \left( P(F) + \lambda^{1+\alpha} \, V_{\alpha}(F) \right).
\end{equation*}
If $\lambda$ is small (the volume of $E$ is small), then the classical perimeter dominates the Riesz potential term, which implies the existence of minimizers. If $\lambda$ is large (the volume of $E$ is large), then the Riesz potential term dominates the classical perimeter, which implies the non-existence of minimizers. 

The plan of the paper is as follows: in Section \ref{sectionClassicalModel}, we review a number of results on the classical liquid drop model, namely, Problem \eqref{classicalLiquidDropModel} with given $m>0$. We first show some results in the case that the kernel $g$ in the energy $\capE_{g}$ is given by the standard one $x \mapsto |x|^{-(N-\alpha)}$ with $\alpha \in (0,\,N)$. Then we show some results in the case of a general kernel $g$. In the general case, we split the section into two parts: the first part is on small volume regime and the second part is on large volume regime. In Section \ref{sectionNonlocalModel}, we review a number of results on the nonlocal extension of Problem \eqref{classicalLiquidDropModel} with given $m>0$.    

\smallskip

{\bf Acknowledgments.} M. Novaga is member of the INDAM-GNAMPA, and was partially supported by the PRIN Project 2019/24.

\section{The classical liquid drop model}\label{sectionClassicalModel}
In this section, we first review some previous works on the classical liquid drop model with the standard kernel $x \mapsto |x|^{-(N-\alpha)}$ of the Riesz potential term. Recently, Frank and Nam in \cite{FrNa} revisited this model and some references are also therein. The main interest from the mathematical point of view is to investigate the following three topics: the existence of minimizer, the non-existence of minimizer, and the minimality of the ball. Kn\"upfer and Muratov in \cite{KnMu01,KnMu02} considered when $g(x) = |x|^{-N+\alpha}$ for $\alpha \in (0,\,N)$ with $N \geq 2$ and proved that there exists constants $0< m_0 \leq m_1 \leq m_2 <\infty$ such that the following three things hold: if $N \geq 2$, $\alpha \in (0,\,N)$, and $m \leq m_1$ , then Problem \eqref{classicalLiquidDropModel} admits a minimizer; if $N \geq 2$, $\alpha \in (N-2,\,N)$, and $m > m_2$, then Problem \eqref{classicalLiquidDropModel} does not admit a minimizer; finally, if $m \leq m_0$, then the ball is the unique minimizer whenever either $N = 2$ and $\alpha \in (N-2,\,N)$, or $3 \leq N \leq 7$ and $\alpha \in (1,\,N)$. Later, Julin in \cite{Julin} proved that, if $N \geq 3$ and $g(x) = |x|^{-(N-2)}$, the ball is the unique minimizer of $\capE_g$ whenever $m$ is sufficiently small. Bonacini and Cristoferi in \cite{BoCr} studied the case of the full parameter range $N \geq 2$ and $\alpha \in (1,\,N)$ when $g(x) = |x|^{-N+\alpha}$. Moreover, for small $\alpha>0$, Bonacini and Cristoferi in \cite{BoCr} gave a complete characterization of the ground state. Namely, they showed that, if $\alpha$ is sufficiently small, there exists a constant $m_c$ such that the ball is the unique minimizer of $\capE_g$ for $m \leq m_c$ and $\capE_g$ does not have minimizers for $m > m_c$. In a slightly different context, Lu and Otto in \cite{LuOt} showed the non-existence of minimizers for large volumes and proved that the ball is the unique minimizer for small volumes when $N = 3$ and $g(x) = |x|^{-1}$. Originally, Lu and Otto were motivated by Thomas-Fermi-Dirac-von Weisz\"aker model (TFDW model) in quantum mechanics and the energy that was studied in \cite{LuOt} includes the background potential, which behaves like an attractive term. After the work by Lu and Otto in \cite{LuOt}, Frank, Nam, and Van Den Bosch in \cite{NaBo,FNB} developed further theory of TFDW model. Finally, we refer to the work by Alama, Bronsard, Choksi, and Topaloglu in \cite{ABCT}, in which they proved that a variant of Gamow's model including the background potential admits minimizers for any volume, due to the effects from the background potential against the Riesz potential. 

In the sequel of this section, we review some results on Problem \eqref{classicalLiquidDropModel} with a more general kernel of the Riesz potential term if the volume is small or large. 


\subsection{Small mass regime with general Riesz kernel}
We now consider the minimization problem with the Riesz potential term associated with a general kernel $g$ under small volume constraint. 

Novaga and Pratelli in \cite{NoPr} actually considered the minimization problem
\begin{equation}\label{minimizationProblemClassicalNovagaPratelli}
	\inf\left\{ \capE^{\varepsilon}_{g}(E) \coloneqq P(E) + \varepsilon \, V_g(E) \mid |E| = |B_1| \right\}
\end{equation}
where $\varepsilon>0$ is a parameter and they showed Theorem \ref{theoremRigidityMinimizersClassicalNovagaPratelli} for sufficiently small $\varepsilon>0$. By a scaling argument, we can observe that Problem \eqref{minimizationProblemClassicalNovagaPratelli} is equivalent to Problem \eqref{classicalLiquidDropModel}. In this setting, Novaga and Pratelli in \cite{NoPr} showed the rigidity of minimizers in two dimension for Problem \eqref{minimizationProblemClassicalNovagaPratelli} for sufficiently small parameter $\varepsilon$. Precisely, they proved the following
\begin{theorem}[Theorem 1.1 in \cite{NoPr}] \label{theoremRigidityMinimizersClassicalNovagaPratelli}
	Assume that $g : \mathR^2 \setminus \{0\} \to (0,\,\infty)$ is a radial, decreasing, positive definite function such that  
	\begin{equation*}
		\int_{B_1} \frac{g(x)}{|x|} \,dx < \infty.
	\end{equation*}
	Then there exists a constant $\varepsilon_0>0$ such that, for any $\varepsilon \in (0,\,\varepsilon_0)$, the balls with volume $|B_1|$ are the unique minimizers for Problem \eqref{classicalLiquidDropModel}.
\end{theorem} 
In addition, they proved the existence of ``generalized'' minimizers of the ``generalized'' energy of $\capE_g$ for any volumes in any dimensions. They precisely study the following minimization problem:
\begin{equation}\label{minimizationProbGeneralizedClassicalEnergy}
	\inf_{H \in \mathN}\inf\left\{ \widetilde{\capE}^{\varepsilon, H}_{g}\left( \{E^k\}_{k} \right) \mid \sum_{k=1}^{\infty} |E^k| = |B_1| \right\}
\end{equation}
where the ``generalized'' energy $\widetilde{\capE}^H_{g}$ is defined as
\begin{equation}
	\widetilde{\capE}^{\varepsilon, H}_{g}\left( \{E^k\}_{k} \right) \coloneqq \sum_{k=1}^{H} \capE^{\varepsilon}_{g}(E_k) 
\end{equation}
for any $\varepsilon > 0$ and $H \in \mathN$. Then one may define a ``generalized'' minimizer of $\widetilde{\capE}^{\varepsilon, H}_{g}$ as a family of sets $\{E^k\}_{k=1}^{H}$ that satisfies
\begin{equation*}
	\widetilde{\capE}^{\varepsilon, H}_{g}\left( \{E^k\}_{k} \right) = \inf_{H \in \mathN}\inf\left\{ \widetilde{\capE}^{\varepsilon, H}_{g}\left( \{E^k\}_{k} \right) \mid \sum_{k=1}^{\infty} |E^k| = |B_1| \right\}
\end{equation*}
and, for each $k \in \{1,\,2, \cdots,\,H\}$, $E^k$ is a minimizer of $\capE^{\varepsilon}_{g}$ among the sets with volume $|E^k|$. Notice that, in the generalized energy, the interaction between different components is not evaluated. This implies that different components can be placed ``at infinity'' from each other. Then Novaga and Pratelli in \cite{NoPr} proved 
\begin{theorem}[Proposition 1.2 in \cite{NoPr}]
	Assume that $g : \mathR^N \setminus \{0\} \to (0,\,\infty)$ is decreasing and satisfies the condition that
	\begin{equation*}
		\int_{B_1}\int_{B_1} g(x-y) \,dx\,dy < \infty.
	\end{equation*}
	Then, for any $\varepsilon > 0$, there exists a generalized minimizer for Problem \eqref{minimizationProbGeneralizedClassicalEnergy}.
\end{theorem}
Theorem \ref{theoremRigidityMinimizersClassicalNovagaPratelli} was improved to the case of higher dimensions by Carazzato, Fusco, and Pratelli in \cite{CFP}. Precisely, they proved the following generalization of Theorem 1:

\begin{theorem}[Theorem A in \cite{CFP}]
	Assume that $g : \mathR^N \setminus \{0\} \to (0,\,\infty)$ is a radial and radially decreasing function such that
	\begin{equation*}
		\int_{B_1} g(x) \,dx < \infty.
	\end{equation*}
	Then there exists a constant $\varepsilon_0>0$ such that, for any $\varepsilon \in (0,\,\varepsilon_0)$, the balls with volume $|B_1|$ are the unique minimizers for Problem \eqref{minimizationProblemClassicalNovagaPratelli}.
\end{theorem}

\subsection{Large mass regime with general Riesz kernel}
Now we consider the minimization problem with the Riesz potential term associated with a general kernel $g$ under large volume constraint. 

Pegon in \cite{Pegon} first showed the existence of minmizers of $\capE_{g}$ under some suitable assumptions on $g$. To state the theorem, we shall introduce the equivalent problem to Problem \eqref{classicalLiquidDropModel} and give some assumptions on $g$. Instead of Problem \eqref{classicalLiquidDropModel}, Pegon in \cite{Pegon} studied the following minimization problem:
\begin{equation}\label{minimizationProblemPegon}
	\inf\left\{ \capE^{\gamma, \lambda}_{g}(E) \coloneqq P(E) - \gamma \, P_{g^1_{\lambda}} (E) \mid |E| = |B_1| \right\}
\end{equation} 
where $\gamma >0$, $g^1_{\lambda}(x) \coloneqq \lambda^{N+1} \, g(\lambda\,x)$ for $\lambda>0$, and $P_{K}$ for a function $K$ is defined by
\begin{equation}\label{generalizedNonlocalPeri}
	P_{K}(E) \coloneqq \int_{E}\int_{E^c} K(x-y)\,dx\,dy
\end{equation}
for a measurable set $E \subset \mathR^N$. Here we use $E^c$ as the complement of $E$. The functional \eqref{generalizedNonlocalPeri} is called the "generalized" nonlocal perimeter associated with the kernel $K$. In this setting, we can indeed show that, whenever $g \in L^1(\mathR^N)$, Problem \eqref{minimizationProblemPegon} is equivalent to Problem \eqref{classicalLiquidDropModel}. Indeed, if $g \in L^1(\mathR^N)$, then we have that
\begin{equation*}
	V_g(E) = |E| \, \|g\|_{L^1(\mathR^N)} - P_g(E)
\end{equation*}
for any $E \subset \mathR^N$ with $|E| < \infty$ and, by a scaling argument, we also have that the volume constraint $|E|=m$ can be replaced with the constraint $|E| = |B_1|$.

To state the results of large volume regime, we now assume that $g$ satisfies the following conditions:
\begin{description}
	\item[(g1)] $g$ is radial, namely, there exists a non-negative function $G: (0,\, \infty) \to \mathR$ such that $g(x) = G(|x|)$.
	
	\item[(g2)] $g$ is integrable and the first moment of $g$ is finite, namely, 
	\begin{equation*}
		\int_{\mathR^N} |x| \, g(x) \,dx < \infty.
	\end{equation*}
\end{description}
With the above assumptions, Pegon in \cite{Pegon} proved
\begin{theorem}[Theorem A in \cite{Pegon}]
	Assume that the kernel $g$ satisfies $(\mathrm{g}1)$ and $(\mathrm{g}2)$ and $\gamma < 1$. Then there exists a constant $\lambda_0>0$ such that, for any $\lambda > \lambda_0$, Problem \eqref{minimizationProblemPegon} admits a minimizer.  
\end{theorem}  
Moreover, Pegon in \cite{Pegon} proved the $\Gamma$-convergence in the $L^1$-topology of $\capE^{\gamma, \lambda}_{g}$ as the volume of minimizers diverges, namely, $\lambda \to \infty$. Precisely, it is shown that 
\begin{equation*}
	\Gamma_{L^1}-\lim_{\lambda \to \infty}\capE^{\gamma, \lambda}(E) = (1 - \gamma)P(E) 
\end{equation*}
for any set $E \subset \mathR^N$ of finite perimeter with $|E| = |B_1|$.

In the case that $g$ is not necessarily in $L^1(\mathR^N)$ and behaves like the kernel of the $s$-fractional perimeter, Mellet and Wu in \cite{MeWu} showed the existence and rigidity of bounded minimizers for small volumes. In addition, They showed that the ball cannot be a global minimizer for large volumes. Precisely, They considered the minimization problem
\begin{equation}\label{minimizationProblemClassicalNonlocalMelletWu}
	\inf\left\{ \capE^{\gamma}_{K}(E) \coloneqq P(E) - \gamma s(1-s) \, P_K (E) \mid |E| = m \right\}
\end{equation}
where $\gamma > 0$, $s \in (0,\,1)$, and $m>0$ under the following assumptions on $K$: 
\begin{description}
	\item[(K1)] $K$ is radially symmetric and non-increasing.
	
	\item[(K2)] $0 \leq K(x) \leq |x|^{-(N+s)}$ for any $x \in \mathR^N$.
\end{description}
Then, Mellet and Wu in \cite{MeWu} first showed the existence of minimizers for small volumes as follows:
\begin{theorem}[Theorem 2.1 in \cite{MeWu}]
	Assume that $K$ satisfies (K1) and (K2). Then there exists a constant $c_0>0$ such that, if $\gamma \, m^{\frac{1-s}{N}} \leq c_0$, then Problem \eqref{minimizationProblemClassicalNonlocalMelletWu} admits a bounded minimizer. 
\end{theorem}
Second, they proved the rigidity of minimizers for small volumes as follows:
\begin{theorem}[Theorem 2.2 in \cite{MeWu}]
	Assume that $K$ satisfies (K1) and (K2) and $2 \leq N \leq 7$. Then there exists a constant $c_1>0$ such that, if $\gamma \, m^{\frac{1-s}{N}} \leq c_1$, then the balls are the unique minimizers for Problem \eqref{minimizationProblemClassicalNonlocalMelletWu}.
\end{theorem}
Finally, with some growth condition of $K$ far away from the origin, they proved that the balls cannot be minimizers for large volumes. Precisely, they further assume the following condition:
\begin{description}
	\item[(K3)] There exist $k_0>0$ and $r_0>0$ such that
	\begin{equation*}
		K(x) \geq \frac{k_0}{|x|^{N+s}} \quad \text{for any $|x| \geq r_0$}.
	\end{equation*}
\end{description}
Then they proved 
\begin{theorem}[Theorem 2.3 in \cite{MeWu}]
	Assume that $K$ satisfies (K1), (K2), and (K3). Then there exists a constant $c_2>0$ such that, if $\gamma \, m^{\frac{1-s}{N}} \geq c_2$, then the balls cannot be minimizers for Problem \eqref{minimizationProblemClassicalNonlocalMelletWu}. 
\end{theorem}

Right after the work by Pegon in \cite{Pegon}, imposing some differentiability of $g$ and some control of its gradient, Merlet and Pegon in \cite{MePe} proved the rigidity of minimizers in two dimension for Problem \eqref{minimizationProblemPegon} for large volumes. They precisely assumed the following condition on $g$:
\begin{description}
	\item[(g3)] $g \in W^{1,1}_{loc}(\mathR^N \setminus \{0\})$ and
	\begin{equation*}
		g(x) = G'(|x|) = \mathcal{O}\left( \frac{1}{|x|^{N+1}} \right)
	\end{equation*}
	as $|x| \to \infty$, where $G$ is as in $(\mathrm{g}1)$ and $\mathcal{O}$ is the big O. 
\end{description}
Then, Merlet and Pegon in \cite{MePe} proved
\begin{theorem}[Theorem A in \cite{MePe}]
	Assume $N=2$ and $g$ satisfies $(\mathrm{g}1)$, $(\mathrm{g}2)$, and $(\mathrm{g}3)$. Then there exists a constant $\lambda_1>0$ such that, for any $\lambda > \lambda_1$, the unit disk is the unique minimizer of Problem \eqref{minimizationProblemPegon}.
\end{theorem}
Without the assumption $(\mathrm{g}3)$, Merlet and Pegon in \cite{MePe} also showed the convexity of minimizers in two dimension as follows:
\begin{theorem}[Theorem 1 in \cite{MePe}] \label{theoremConvexityMiniProbMerletPegon}
	Assume $N=2$ and $g$ satisfies $(\mathrm{g}1)$ and $(\mathrm{g}2)$. Then there exists a constant $\lambda_2>0$ such that, for any $\lambda > \lambda_2$, every minimizer for Problem \eqref{minimizationProblemPegon} is convex.
\end{theorem}
The idea of the proof of Theorem \ref{theoremConvexityMiniProbMerletPegon} is to properly slice the minimizer of $\capE_{g}$ with respect to some line and reduce the argument to the study of the critical energy on the real lines. This strategy is specific to the two dimensional case and it is not applicable to the case of higher dimensions. The convexity of minimizers with large volumes is still open in higher dimensions. However, Merlet and Pegon in \cite{MePe} showed that, in any dimension and for sufficiently large $\lambda>0$, the balls are the unique minimizers for Problem \eqref{minimizationProblemPegon} among the sets that are ``close'' to balls in a proper sense.

\section{The nonlocal liquid drop model}\label{sectionNonlocalModel}
In this section, we study a nonlocal version of the classical liquid drop model and we see how different the classical and nonlocal problems are. The nonlocal model that we will investigate aims at dealing with long-range interactions between particles in physics. The long-range interactions can influence the existence of minimizers and, indeed, this contribution from ``far-away'' can produce some distinct phenomena from the classical problems.   

Let $m > 0$ and $s \in (0,\,1)$. We consider the minimization problem
\begin{equation}\label{nonlocalLiquidDropModel}	
	\inf\left\{ \capE_{s,g}(E) \coloneqq P_s(E) + V_g(E) \mid |E| = m \right\}
\end{equation}
where $P_s(E)$ is the $s$-fractional perimeter of a set $E \subset \mathR^N$ and $V_g(E)$ is the Riesz potential term of a set $E \subset \mathR^N$ associated with the kernel $g$ (see \eqref{definitionGeneralRieszPotential} for the definition). The definition of the $s$-fractional perimeter is given as
\begin{equation}\label{defNonlocalPerimeter}
	P_s(E) \coloneqq \int_{E}\int_{E^c} \frac{1}{|x-y|^{N+s}} \,dx\,dy
\end{equation}
for any measurable set $E \subset \mathR^N$. The notion of the $s$-fractional perimeter was introduced by Caffarelli, Roquejoffre, and Savin in \cite{CRS} to study the classical phase-field model with long-range correlations.

One can easily see that Problem \eqref{nonlocalLiquidDropModel} can be regarded as a nonlocal version of Problem \eqref{classicalLiquidDropModel} and, from the fact that
\begin{equation*}
	(1-s) \, P_s(E) \xrightarrow[s \uparrow 1]{} |\partial B_1| \, P(E)
\end{equation*}
for any smooth set $E \subset \mathR^N$ (see \cite{ADPM}), Problem \eqref{nonlocalLiquidDropModel} somehow approximates to Problem \eqref{classicalLiquidDropModel} in $s \in (0,\,1)$. As is discussed in the previous section, by a scaling argument, one can heuristically observe the existence of minimizers for small volumes and the non-existence of minimizers for large volumes. In the sequel, we shall see that the heuristic argument is indeed valid in Problem \eqref{nonlocalLiquidDropModel} and review a number of results on the nonlocal liquid drop model with the standard or general kernel of the Riesz potential term. 

\subsection{Small mass regime}
We here present the existence and rigidity results of minimizers with small volumes in the case that the kernel of the Riesz potential term is given by the standard one. Precisely, we consider Problem \eqref{nonlocalLiquidDropModel}, provided that the kernel $g$ is given by $x \mapsto |x|^{-(N-\alpha)}$ with $\alpha \in (0,\,N)$.

Figalli, Fusco, Maggi, Millot, and Morini in \cite{FFMMM} firstly investigated Problem \eqref{nonlocalLiquidDropModel} when $g(x) = |x|^{-(N-\alpha)}$ with $\alpha \in (0,\,N)$ and $m$ is sufficiently small. They proved 
\begin{theorem}[Theorem 1.3 in \cite{FFMMM}] \label{resultbyFFMMM}
	Let $N \geq 2$ and $s \in (0,\,1)$. Assume that $g(x) = |x|^{-(N-\alpha)}$ with $\alpha \in (0,\,N)$. Then there exists a constant $m_* > 0$ such that, for any $m \in (0,\,m_*)$, Problem \eqref{nonlocalLiquidDropModel} admits the balls as the unique minimizers up to translations.
\end{theorem}     
Moreover, Figalli, et al. in \cite{FFMMM} also gave the threshold of volume-constrained $L^1$-local minimality of balls for $\capE_{s,g}$. See \cite[Theorem 1.5]{FFMMM} for the precise claim.

If the kernel $g$ of the Riesz potential term is given in a general way, Carazzato in \cite{Carazzato} showed that the balls are the unique minimizers for Problem \eqref{nonlocalLiquidDropModel} for small volumes. Carazzato in \cite{Carazzato} precisely studied the following minimization problem:
\begin{equation}\label{minimizationProblemNonlocalCarazzato}
	\inf\left\{ \capE^{\varepsilon}_{s,g}(E) \coloneqq P_s(E) + \varepsilon \, V_g(E) \mid |E| = |B_1| \right\}
\end{equation}
where $\varepsilon >0$ is a parameter. Notice that Problem \eqref{minimizationProblemNonlocalCarazzato} is a nonlocal version of Problem \eqref{minimizationProblemClassicalNovagaPratelli}. In this setting, the following assumptions are considered:
\begin{description}
	\item[(H1)] $g$ is non-negative and radially non-increasing, namely,
	 \begin{equation*}
	 	g(\lambda \, x) \leq g(x) \quad \text{for any $x \in \mathR^N \setminus \{0\}$ and $\lambda \geq 1$.}
	 \end{equation*}
 
 	\item[(H2)] $g \in L^1_{loc}(\mathR^N)$ and there exists a constant $R_g>0$ such that $g$ is bounded in $\{|x| > R_g\}$. 
\end{description}
Note that the standard kernel $|x|^{-(N-\alpha)}$ with $\alpha \in (0,\,N)$ obviously satisfies (H1) and (H2). Then Carazzato in \cite{Carazzato} proved
\begin{theorem}[Theorem A in \cite{Carazzato}]
	Let $s \in (0,\,1)$. Assume that $g$ satisfies (H1) and (H2). Then there exists a constant $\varepsilon_0 > 0$ such that, for any $\varepsilon < \varepsilon_0$, the balls are the unique minimizers for Problem \eqref{minimizationProblemNonlocalCarazzato}. 
\end{theorem}

\subsection{Large mass regime}
As an analogy of the classical problem, the second author of the present article showed in \cite{Onoue} that there exists no minimizer of $\capE_{s,g}$ for large volumes whenever $N \geq 2$ and $g(x)=|x|^{-1}$. The author studied $\capE_{s,g}$ in more general settings with the background potential. Precisely, given $\mu \geq 0$ and $\beta \in [0,\,N+1)$, the author studied the minimization problem
\begin{equation}\label{nonlocalEnergywithBackgroundPotential}
	\inf\left\{ \capE_{K, g, \mu, \beta}(E) \coloneqq P_K(E) + V_g(E) - \mu\int_{E}\frac{dx}{|x|^{\beta}} \mid |E| = m \right\}
\end{equation}   
where $K : \mathR^N \setminus \{0\} \to [0,\,\infty)$ is a measurable function which, roughly speaking, behaves as a function $x \mapsto |x|^{-(N+s)}$ with some $s \in (0,\,1)$ and $P_{K}(E)$ of $E \subset \mathR^N$ is defined in \eqref{generalizedNonlocalPeri}. 

Notice that, if $K(x) = |x|^{-(N+s)}$ and $\mu = 0$, then Problem \eqref{nonlocalEnergywithBackgroundPotential} coincides with Problem \eqref{nonlocalLiquidDropModel}. In this general setting, the author in \cite{Onoue} proved
\begin{theorem}[Theorem 2.2 in \cite{Onoue}]\label{theoremNonexistenceMiniLargeMassNonlocal}
	Let $N \geq 2$, $s \in (0,\,1)$, and $\beta=1$ and let $g(x) = |x|^{-1}$. Assume that $K$ is radially symmetric and suitably controlled by $x \mapsto |x|^{-(N+s)}$ with $s \in (0,\,1)$. Then there exists a constant $m_0 > 0$ (explicitly given) such that, for any $m \geq m_0$, Problem \eqref{nonlocalEnergywithBackgroundPotential} admits no bounded minimizer.
\end{theorem} 
Note that, if we further impose on the kernel $K$ some control of oscillation, then we can show that Problem \eqref{nonlocalEnergywithBackgroundPotential} admits no minimizer for large volumes. See \cite[Corollary 2.3]{Onoue}.
   
Problem \eqref{nonlocalEnergywithBackgroundPotential} can be also regarded as a nonlocal generalization of the work by Lu and Otto shown in \cite{LuOt,LuOt02}. Lu and Otto were motivated by a sharp interface version of Thomas-Fermi-Dirac-von Weisz\"aker model (TFDW model) in quantum mechanics. The problem is also related to the ionization conjecture, which states that the number of electrons that can be bound to an atomic nucleus of charge $\mu$ cannot exceed $\mu+1$. For the readers interested in  the TFDW model or ionization conjecture, we refer to the non-exhaustive list of references \cite{BBL,BeLi,Lieb02,LuOt,LSST,Nam,NaBo,Sigal,Solovej,Solovej02}. 

On the contrary, the existence of minimizers for Problem \eqref{nonlocalLiquidDropModel} for large volumes has been studied by the authors of this article under a ``good'' control of the kernel $g$. One may heuristically observe that, if the Riesz potential term of $\capE_{s,g}$ is controlled in a good way by the $s$-fractional perimeter term, then there may be the possibilities that minimizers exist for large volumes. Indeed, we can rigorously show that this heuristic argument is valid. To see this, we assume that the kernel $g$ of the Riesz potential term satisfies the following conditions:
\begin{description}
	\item[(g1)'] $g$ is non-negative and radially non-increasing, namely,
	\begin{equation*}
		g(\lambda \, x) \leq g(x) \quad \text{for any $x \in \mathR^N \setminus \{0\}$ and $\lambda \geq 1$.}
	\end{equation*}
	
	\item[(g2)'] $g$ is symmetric at the origin, namely, $g(-x) = g(x)$ for any $x \in \mathR^N$.
	
	\item[(g3)'] There exists $R>0$ and $\eta \in (0,\,1)$ such that
	\begin{equation*}
		g(x) \leq \frac{\eta}{|x|^{N+s}} \quad \text{for any $|x| \geq R$.}
	\end{equation*} 
\end{description}
Note that $(\mathrm{g}1)'$ is same as Assumption (H1) in the previous subsection. Given these assumptions, the authors of this article proved the following existence result:
\begin{theorem}[Theorem 2.3 in \cite{NoOn}]\label{theoremExistenceMiniAnyVolumes}
	Assume that $g$ satisfies $(\mathrm{g}1)'$, $(\mathrm{g}2)'$, and $(\mathrm{g}3)'$. Then Problem \eqref{nonlocalLiquidDropModel} admits minimizers for any volume $m>0$.
\end{theorem}  

The idea of the proof is based on the concentration-compactness lemma introduced by P.L. Lions in \cite{Lions01,Lions02}. In general, when one has a minimizing sequence of Problem \eqref{nonlocalLiquidDropModel}, one has three possibilities: the first one is that the sets of the minimizing sequence tend to be the union of small ``dusts'' and vanish in the limit. This is so-called {\it Vanishing} in the sense of Lions. The second one is that the sets of the sequence tend to split into pieces, which are placed far away from each other, in the limit. 
This is the so-called {\it Dichotomy} in the sense of Lions. The last possibility is that the sequence converges to some set, which turns out to be a minimizer of the energy. This is so-called {\it Compactness} in the sense of Lions. Basically in our problem, we can exclude the possibility of {\it Vanishing} thanks to the uniform bound and isoperimetric inequality of the $s$-fractional perimeter. In addition, thanks to the ``good'' control of the kernel $g$, we can also exclude the possibility of {\it Dichotomy}. Therefore, we can obtain {\it Compactness} in the sense of Lions and this implies the existence of minimizers.

As is mentioned above, one may not exclude the possibility of {\it Dichotomy} without some proper control of the kernel $g$ and thus one does not have the existence of minimizers. However, one may obtain the so-called ``generalized'' minimizer of a ``generalized'' energy of $\capE_{s,g}$. Before stating our result, we introduce the notion of a generalized minimizer and a generalized energy. Note that, in Section \ref{sectionClassicalModel}, the similar notion is given for the classical liquid drop model. For any $m>0$. we define the ``generalized'' energy of $\capE_{s,g}$ over a family of sequences of sets $\{E^k\}_{k\in\mathN}$ with $\sum_{k}|E^k| = m$ as
\begin{equation}\label{generalizedNonlocalEnergy}
	\widetilde{\capE}_{s,g}\left( \{E^k\}_{k} \right) \coloneqq \sum_{k=1}^{\infty} \capE_{s,g}(E^k).
\end{equation}
With this energy, we define the ``generalized'' minimizers by a family of sets that solves the minimization problem
\begin{equation}\label{minimizationProbGeneralizedNonlocalEnergy}
	\inf\left\{ \widetilde{\capE}_{s,g}\left( \{E^k\}_{k} \right) \mid \sum_{k=1}^{\infty} |E^k|=m \right\}
\end{equation}
and, for each $k \in \mathN$, $E^k$ is a minimizer of Problem \eqref{nonlocalLiquidDropModel} with $m = |E^k|$. By definition, one may have the possibility that the minimizing sequence of Problem \eqref{minimizationProbGeneralizedNonlocalEnergy} tends to split into pieces, which are placed at an infinite distance from each other. Thus we may not be able to exclude {\it Dichotomy}. 

To see the existence of generalized minimizers, we assume, instead of $(\mathrm{g}3)'$, the following condition on $g$, which is much weaker than $(\mathrm{g}3)'$:
\begin{description}
	\item[(g4)'] $g$ vanishes at infinity, namely, $g(x) \to 0$ as $|x| \to \infty$. 
\end{description}
With Assumption $(\mathrm{g}4)'$, we can prove the existence of the generalized minimizers of $\widetilde{\capE}_{s,g}$ for Problem \eqref{minimizationProbGeneralizedNonlocalEnergy} as follows:
\begin{theorem}[Theorem 2.4 in \cite{NoOn}] \label{existenceGeneralizedMiniGeneralizedNonlocalEnergy}
	Assume that $g$ satisfies $(\mathrm{g}1)'$, $(\mathrm{g}2)'$, and $(\mathrm{g}4)'$. Then Problem \eqref{minimizationProbGeneralizedNonlocalEnergy} admits generalized minimizers for any volume $m>0$. 
\end{theorem}

We now consider the asymptotic behavior of minimizers as the volume goes to infinity under the ``good'' control of $g$. To see this, we shall play with another minimization problem which is equivalent to Problem \eqref{nonlocalLiquidDropModel}. Given $\lambda>0$, we set the rescaled kernel $g^s_{\lambda}$ as $g^s_{\lambda}(x) \coloneqq \lambda^{N+s}  g(\lambda\,x)$ for any $x \in \mathR^N \setminus \{0\}$. Then we study the minimization problem
\begin{equation}\label{equivalentMinimizationProbNonlocalEnergy}
	\inf\left\{ \widehat{\capE}^{\lambda}_{s,g}(E) \coloneqq P_s(E) - P_{g^s_{\lambda}}(E) \mid  |E| = |B_1| \right\}
\end{equation}
where $P_{g^s_{\lambda}}$ is given in \eqref{generalizedNonlocalPeri} with $K = g^s_{\lambda}$ and we see how the minimizers of $\capE^{\lambda}_{s,g}$ behave as $\lambda \to \infty$. From the same argument shown in Section \ref{sectionClassicalModel}, one can easily observe that Problem \eqref{equivalentMinimizationProbNonlocalEnergy} is equivalent to Problem \eqref{nonlocalLiquidDropModel} whenever $g \in L^1(\mathR^N)$. Under this situation, we can show that the sequence of solutions $\{E_{\lambda}\}_{\lambda}$ to Problem \eqref{equivalentMinimizationProbNonlocalEnergy} converges to the Euclidean ball with volume $|B_1|$ as $\lambda \to \infty$. To see this, we assume that the kernel $g$ satisfies the following condition, which is much stronger than $(\mathrm{g}3)$:
\begin{description}
	\item[(g5)'] There exists $\eta \in (0,\,1)$ such that 
	\begin{equation*}
		g(x) \leq \frac{\eta}{|x|^{N+s}} \quad \text{for $x \in \mathR^N \setminus \{0\}$}, \quad g(x) = o\left( \frac{1}{|x|^{N+s}} \right) \quad \text{as $x \to \infty$}. 
	\end{equation*}
\end{description}
Notice that $(\mathrm{g}5)'$ does not necessarily require the kernel $g$ to be integrable in $\mathR^N$. Taking the above into account, we can obtain    
\begin{theorem}[Theorem 2.6 in \cite{NoOn}]\label{theoremAsymptoticMinimizersEquivalentProb}
	Let $\{\lambda_i\}_{i \in \mathN} \subset (0, \, \infty)$ be such that $\lambda_i \to \infty$ as $i \to \infty$ and let $\{F_i\}_{i \in \mathN}$ be a sequence of minimizers of $\widehat{\capE}^{\lambda_i}_{s,g}$ with $|F_i| = |B_1|$. Assume that $g$ satisfies $(\mathrm{g}1)'$, $(\mathrm{g}2)'$, and $(\mathrm{g}5)'$. Then, up to translations,  
	\begin{equation*}
		|F_i \Delta B_1| \xrightarrow[i \to \infty]{} 0. 
	\end{equation*} 
\end{theorem}
Moreover, we can show that $\widehat{\capE}^{\lambda}_{s,g}$ converges, in the sense of $\Gamma$-convergence, to the $s$-fractional perimeter as $\lambda \to \infty$. 

As concluding remarks of the article, we present several open problems on the nonlocal extension of the classical liquid drop model. As pursued in the classical liquid drop model, one can come up with the following questions: first, one may ask
\begin{description}
	\item[Q1] Does Problem \eqref{nonlocalLiquidDropModel} admit no minimizers for large volumes, if the kernel $g$ is the standard one, namely, $g(x) = |x|^{-(N-\alpha)}$ with $\alpha \in (0,\,N)$? 
\end{description}
As shown in Theorem \ref{theoremNonexistenceMiniLargeMassNonlocal}, we have the non-existence result only in the case of the specific kernel $g$. Second, one may ask
\begin{description}
	\item[Q2] What is the critical volume after which Problem \eqref{nonlocalLiquidDropModel} does not admit minimizers (or equivalently the critical volume below which Problem \eqref{nonlocalLiquidDropModel} admits minimizers)?
\end{description}
This question is still open even in the case of the classical liquid drop model. See \cite{FrNa} for the detail. Finally, one may ask
\begin{description}
	\item[Q3] Are balls the unique solutions to Problem \eqref{nonlocalLiquidDropModel} under suitable assumptions on $g$ for large volumes?
\end{description}
In the local case, such rigidity result was proved in \cite{MePe} in two dimensions.
As shown in Theorem \ref{theoremAsymptoticMinimizersEquivalentProb}, we only have that the minimizers for Problem \eqref{equivalentMinimizationProbNonlocalEnergy} (or equivalently Problem \eqref{nonlocalLiquidDropModel} if $g \in L^1$) converge to the balls as the volume diverges in the $L^1$-topology.


%
%

\end{document}